\documentclass[12pt]{article}
\usepackage{amsmath,amssymb, amsfonts, amsthm,amscd}
\usepackage[T2A]{fontenc}
\usepackage[cp1251]{inputenc}
\usepackage{graphicx}

\theoremstyle{plain}
\newtheorem{thm}{Theorem}
\newtheorem{prop}{Proposition}

\theoremstyle{definition}
\newtheorem{defini}{Definition}
\newtheorem{defin}[defini]{Definition}

\newenvironment{Proof}
    {\par\noindent {\bf Proof.} \quad}{\par\hfill $\Box$}
\newcommand{\Ann}{\mathrm{Ann}_r}
\newcommand{\ov}{\overline}
\newcommand{\ovR}{\overline{R}}
\newcommand{\ovb}{\overline{b}}
\newcommand{\ova}{\overline{a}}
\newcommand{\ovc}{\overline{c}}
\newcommand{\ovd}{\overline{d}}
\newcommand{\ovs}{\overline{s}}
\newcommand{\ovr}{\overline{r}}
\newcommand{\ovv}{\overline{v}}
\newcommand{\ovu}{\overline{u}}
\newcommand{\ove}{\overline{e}}
\newcommand{\ovx}{\overline{x}}
\newcommand{\la}{\lambda}
\newcommand{\al}{\alpha}
\newcommand{\be}{\beta}
\newcommand{\ovla}{\overline{\lambda}}
\newcommand{\ovbe}{\overline{\beta}}
\newcommand{\ovmu}{\overline{\mu}}

\renewcommand\refname{\large \textbf{References}}

\begin{document}
\begin{center}
{\huge Rings of Dyadic range 1}
\end{center}
\vskip 0.1cm \centerline{{\Large Bohdan Zabavsky}}

\vskip 0.3cm

\centerline{\footnotesize{Department of Mechanics and Mathematics,  Ivan Franko National University of Lviv,  Ukraine}}
\vskip 0.5cm

\centerline{\footnotesize{February, 2017}}
\vskip 0.7cm

\footnotesize{\noindent\textbf{Abstract:} \textit{
      Using the concept of ring  diadic range 1 we proved that a commutative Bezout ring is an elementary divisor ring iff it is a ring diadic range 1.} }


\vspace{1ex}
\footnotesize{\noindent\textbf{Key words and phrases:} \textit{elementary divisor ring, diadic range 1, stable range 1, stable range 2.}

}

\vspace{1ex}
\noindent{\textbf{Mathematics Subject Classification}}: 06F20,13F99.

\hskip 1,5 true cm

\normalsize

\section{Introduction}
\label{g1220}

One of the main sources of almost all researches in the present paper is a rather problem of a full description  of the elementary divisor rings. The notion of an elementary divisor ring was introduced by Kaplansky \cite{KaplEDMod}. Recall that a matrix over an associative ring has a canonical diagonal reduction if can be reduced to a diagonal form by left and right multiplication by some invertible matrices of the corresponding sizes and so that each diagonal element is a full divisor of the following one. If any matrix over a ring has a canonical diagonal reduction then such a ring is called an elementary divisor ring \cite{KaplEDMod}. In the same paper, Kaplansky showed that any finitely presented module on an elementary divisor ring can be decomposed into a direct sum of cyclic modules. In the case of commutative rings the reverse statement is proved. Namely, if any finitely presented module over a ring can be decomposed into a direct sum of cyclic modules then this ring is an elementary divisor ring \cite{LarsLevShorEDR}. This result is a partial solution to the problem of Warfield \cite{Warf}.

There are a lot of researches that deal with the matrix diagonalization in different casses (the most comprehensive history of these researches can be found in \cite{Zabavsk}).

Specific role in modern research on elementary divisor rings is played by a K-theoretical invariant such as the stable range that was established in 1960 by Bass \cite{BASS64}. One of the most fruitful aspects of Bass studies was the following fact: a lot of answers to the problems of the linear algebra over rings becomes simpler if we increase the dimension of the considered object (the rank of the projective module, the size of the matrix etc.) and, furthermore, the answer is independent of the choice of the base ring for rather big dimensions values, as well as is independent of the dimension of the current object -- in functions only in terms of the "geometry" of a given module. Moreover, it has been discovered that in the commutative ring case there exist structural theorems on the these objects starting from some small values of a stable range (for example, 1 or 2) which depends only on the considered problem, but does not depend on neither the dimension or structure of the base ring.

For example, an important role in studying of the elementary divisor rings is played by the Hermitian rings. A ring is called right (left) Hermitian if all $1 \times 2$ ($2 \times 1$) matrices over this ring have diagonal reduction over this ring. An Hermitian ring is a ring which is both right and left Hermitian \cite{KaplEDMod}.

Note that any Hermitian ring is a finitely generated principal ideal ring (will continue to call them Bezout ring) i.e. a ring in which any finitely generated left or right ideal is principal \cite{KaplEDMod}. In the case of the commutative rings we have the following result.

\begin{thm}\cite{Zabavsk}\label{t1}
      A commutative Bezout ring is an Hermitian ring if and only if it is a ring of stable range 2.
\end{thm}

\begin{defin}
      We say that a ring $R$ has stable range 2 if for any elements $a, b, c \in R$ the equality $aR + bR + cR = R$ implies that there are some elements $\la, \mu \in R$ such that
      $(a + c\la)R + (b + c\mu)R = R\cite{Zabavsk}$.
\end{defin}

Under a ring in this paper we always understand an associative ring with unit and $1 \neq 0$.

\begin{defin}
      A commutative ring $R$ is called a Gelfand ring if for any $a, b \in R$ such that $a + b = 1$ there exist elements $x, y$ such that $(1 + ax)(1 + by) = 0$ \cite{Cont_PMR}.
\end{defin}

\begin{defin}
      Let $R$ be a commutative Bezout domain. We say that $R$ is a ring of Gelfand range 1 if for any $a, b \in R$ such that $aR + bR = R$ there exists an element $\la \in R$ such that the factor-ring $R/(a + b\la)R$ is a Gelfand ring \cite{Zab2017}.
\end{defin}

The main result is thefollowing theorem.

\begin{thm}\cite{Zab2017}\label{t2}
      A commutative Bezout domain is an elementary divisor ring if and only if it is a ring of Gelfand range 1.
\end{thm}

This result give a solution of the problem of elementary divisor rings for different classes of a commutative Bezout domains, in particular, for a PM$^\star$ ring local Gelfand domains an so on \cite{Pihura_KMP}.

In this paper based on the concept of a ring dyadic range 1 we similarly describe commutative elementary divisor rings.

\section{Ring dyadic range 1}

Let $R$ be an assotiative ring with unit and $1 \neq 0$.

\begin{defin}
      Let $a, b \in R$ and $aR + bR = R$. We say that a pair $(a, b)$ has a right diadem (or pair $(a, b)$ is right dyadic) if there is an element $\la \in R$ such that for the element $a + b\la$ and any elements $c, d \in R$ such that $(a + b\la)R + cR + dR = R$ there is an element $\mu \in R$ such that $(a + b\la)R + (c + d\mu)R = R$.
      We call the element $a + b\la$ a right diadem of the pair $(a, b)$.
      By analogy, one can introduce the notion of left diadem and left diadic pair. A right diadem which is a left diadem will called a diadem.
\end{defin}

An obvious example of a diadic pair is the so called trivial diadic pair $(a, u)$, where $u$ is an invertible element of a ring $R$ and $a$ is any element of $R$, $u + a0$ and $a + (-au^{-1} + 1)u$ are right diadems of the pair $(a, u)$.

A nontrivial example is a pair $(a, a + u)$ where $a \in R$ and $u$ is an invertible element of $R$. Here $a + (a + u) - 1$, $(a + u) + a(-1)$ are right diadems of the pair $(a, a + u)$.

\begin{defin}
      If for a right diadic pair $(a, b)$ there exists an element $\la \in R$ such that $a + b\la$ is an invertible element, then the diadem $a + b\la$ is called a trivial diadem.
\end{defin}

\begin{defin}
      We say that a ring $R$ has stable range 1 if for any elements $a, b \in R$ the equality $aR + bR = R$ implies that there is some element $\la \in R$ such that $(a + b\la)R = R$.
      If such element $\la \in R$ always can be taken to be idempotent then we say that $R$ has the idempotent stable range 1.
\end{defin}

It is well known that every semiperfect ring, unit regular ring and strongly $\pi$-regular ring has stable range 1. Meanwhile, every left (right) quasi-duo exchange ring has stable range 1 and every exchange ring of bounded index of nilpotency has stable range 1.

Chen \cite{Ch9} showed that an Abelian ring is clean if and only if it has idempotent stable range 1. Note that a ring of stable range 1 any right dyadic pair $(a, b)$ has a trivial diadem.

\begin{defin}
      We say that a ring $R$ is a ring of right dyadic range 1 if for any elements $a, b$ the equality $aR + bR = R$ implies that the pair $(a, b)$ has a right diadem. Similarly, we define determined a ring of a left dyadic range 1. A ring of right dyadic range 1 which is a ring of left dyadic range 1 is called a ring of dyadic range 1.
\end{defin}

An obvious example of a ring of dyadic range 1 is a ring of stable range 1. Moreover, the following result holds.

\begin{thm}\label{t3}
      Any Bezout ring of right dyadic range 1 is a ring of stable range 2.
\end{thm}

\begin{Proof}
      Let $aR + bR + cR = R$. Since $R$ is a right Bezout ring, $bR + cR = dR$ and $aR + dR = R$. Let $\nu = a + d\la$ be a right diadem of pair $(a, d)$. Then $\nu = a + bx + cy$ for some elements $x, y \in R$. Note that $\nu R + bR + cR = R$. Since $aR + bR + cR = R$. According to the definition of ring $R$ of right dyadic range 1, we have $\nu R + (b + c\mu)R = R$ for some element $\mu \in R$. Then $\nu s + (b + c\mu)t = 1$ for some elements $s, t \in R$. Note, that $Rs + Rt = R$. Since $\nu = a + bx + cy$, we see that $(a + bx + cy)s + (b + c\mu)t = 1$ and $as + b(xs + t) +c(ys + \mu t) = 1$. Since $Rs + Rt = R$, we obtain $Rs + R(xs + t) = R$. Indeed, if $Rs + R(xs + t) = Rh$, where $h$ is not an invertible element of $R$, then $s = s_0h$, $xs + t = \la h$ for some elements so $\la \in R$. Then $\la h = xs + t - xs_0h + t$ and $t = (\la - xs_0)h$, which is impossible, since $Rs + Rt = R$. Thus $Rs + R(xs + t) = R$ whence $us + \mu (xs + t) = ys + \mu t$ for some elements $u, v \in R$.
      Since $as + b(xs + t) + c(ys + \mu t) = 1$ and $us + \nu (xs + t) = ys + \mu t$, we have $as + b(xs + t) + cus + c\nu(xs + t) = 1$ that $(a + cu)s + (b + c\nu)(xs + t) = 1$ and $(a + cu)R + (b + c\nu)R = R$. Thus, $R$ is a ring of stable range 2.
\end{Proof}

\begin{prop}\label{p4}
      A ring $R$ is a ring of right dyadic range 1 if and only if $R/J(R)$ is a ring of right dyadic range 1, where $J(R)$ is the Jacobson radical of $R$.
\end{prop}

\begin{Proof}
      Let $\ovR = R/J(R)$ and $\ova\ovR + \ovb\ovR = \ovR$, where $\ova = a + J(R)$, $\ovb = b + J(R)$. Since $\ova\ovR + \ovb\ovR = \ovR$, we obtain $aR + bR = R$. Since $R$ is a ring of right dyadic range 1 then there exists an element $\la \in R$ such that $a + b\la$ is a right diadem pair $(a, b)$. Let $\ovla = \la + J(R)$ and $(\ova + \ovb\ovla)\ovR + \ovc\ovR + \ovd\ovR = \ovR$ where $\ovc = c + J(R)$, $\ovd = d + J(R)$. Then $(a + b\la)R + cR + dR = R$ and since $a + b\la$ is right diadem pair $(a, b)$. We have $(a + b\la)R + (c + d\mu)R = R$ for some element $\mu \in R$ and $(\ova + \ovb\ovla)\ovR + (\ovc + \ovd\ovmu)\ovR = \ovR$, thus $\ova + \ovb\ovla$ is a right diadem pair $(\ova, \ovb)$. If $\ovR$ is right of right dyadic range 1, then the fact that $\ovR = R/J(R)$ is a ring of right diadic range 1 is obviously.
\end{Proof}

\begin{prop}\label{p5}
      Let $R$ be a ring of right diadic range 1 and $aR = bR$. Then there exist diadems $d_1, d_2$ for some dyadic pairs of $R$ such that $ad_1 = b$, $bd_2 = a$.
\end{prop}

\begin{Proof}
      Since $aR = bR$, we have $a = bs$, $b = at$ for some elements $s, t \in R$. This implies $a(1 - ts) = 0$ and $1 - ts \in \Ann(a)$. Note that if $\Ann(a) = (0)$, then $ts = 1$. By Theorem \ref{t3}, $R$ is a ring of stable range 2, by \cite{Ch9} $R$ is finite, i.e. $t, s$ is invertible ring and proposition is proved.

      Let $\Ann(a) \neq (0)$. Since $1 - ts \in \Ann(a)$, then $tR + \Ann(a) = R$. From here, we have $tx + \la = 1$ for some elements $x \in R$ and $x \in \Ann(a)$. Since $R$ is a ring of right diadic range 1, let $t + \al\la = d_1$ is diadem pair $(t, \al)$. Then $at + a\al\la = ad_1$, since $a\al = 0$, and $at = b$, we have $b = ad_1$. Which have similar $bd_2 = a$ for some right diadem.
\end{Proof}

Moreover, for a class right quasi morphic ring we have the inverse statement.

\begin{prop}\label{p6}
      Let $R$ be a right quasi morphic ring in which the condition $aR = bR$ followed that $a = bd_1$, $b = ad_2$ for some diadems $d_1, d_2 \in R$. Then $R$ is a ring of right dyadic range 1.
\end{prop}

\begin{Proof}
      Let $xR + yR = 1$, then $xz -1 \in yR$ for some element $z \in R$. Let $yR = \Ann(\al)$ and $\al xR = \Ann(\be)$ for some elements $\al, \be \in R$. The following elements are since $R$ is right quasi morphic ring. Note that for any element $r \in R$ we have $(\be\al)xr = \be(\al x)r = 0$ and $(\be\al)yr = \be(\al y)r = 0$ since $\al(yr) = 0$ for any element $r \in R$. So $xR \subseteq \Ann(\be\al)$ and $yR \subseteq \Ann(\be\al)$. Since $xR + yR = R$, we have $1 \in \Ann(\be\al)$, i.e. $\be\al = 0$. So $\al \in \Ann(\be)$, i.e. $\al R \subseteq \Ann(\be)$. Also we have $\Ann(\be) = \al xR \subset \al R$.

      Therefore we have $\Ann(\be) = \al xR = aR$. Under the conditions imposed on $R$, we have with conditions $dxR = \al R$ follows $\al x = \al d$ for some right diadem $d \in R$. So $d(x - d) = 0$ i.e. $x - d \in \Ann(a) = yR$. This implies $x + y\la = d$ for some elements $\al \in R$ that $R$ is ring of right dyadic range 1.
\end{Proof}

The following proposition will be useful in the sequel.

\begin{prop}\label{p7}
      Let $R$ be a commutative ring and a pair $(a, b)$ be a dyadic pair. The element $a + b\la$ is diadem if and only if factor-ring $R/(a + b\la)R$ is a ring of stable range 1.
\end{prop}

\begin{Proof}
      Let $a + b\la$ be a diadem. Denote $\ovR = R/(a + b\la)R$ and let $\ovc\ovR + \ovd\ovR = \ovR$, where $\ovc = c + (a + b\la)R$, $\ovd = d + (a + b\la)R$. Since $\ovc\ovR + \ovd\ovR = \ovR$ then $(a + b\la)R + cR + dR = R$. Since $a + b\la$ is a diadem, there exist $\mu \in R$ such that $(a + b\la)R + (c + d\mu)R = R$. Whence $(\ovc + \ovd\ovmu)\ovR = \ovR$, where $\ovmu = \mu + (a + b\la)R$. Thus, it is proved that the stable range of $\ovR = R/(a + b\la)R$ is equal to 1. Sufficiency is obvious.
\end{Proof}

\begin{prop}\label{p8}
      Let $R$ be a commutative ring and let $c \in R\setminus \{0\}$. If for any $a, b \in R$ such that $aR + bR + cR = R$ there exist $r, s \in R$ such that $c = rs$ and $rR + sR = R$, $rR + aR = R$, $sR + bR = R$ then $R/cR$ is an exchange ring.
\end{prop}

\begin{Proof}
      Denote $\ovR = R/cR$. Since $aR + bR + cR = R$, then $\ovb\ovR + \ova\ovR = \ovR$ where $\ovb = b + cR$, $\ova = a + cR$. Denote $\ovr = r + cR$, $\ovs = s + cR$. Since $rR + sR = R$, one has $ru + sv = 1$. We have $\ovr^2\ovu = \ovr\ovu$, $\ovs^2\ovv = \ovs\ovv$, where $\ovu = u +cR$, $\ovv = v +cR$. Let $\ovs\ovv = \ove$, obviously $\ove^2 = \ove$ and $\ov1 - \ove = \ovr\ovu$. Since $rR + aR = R$, we obtain $r\al + a\be = 1$ for some elements $\al, \be \in R$, then $rsv\al + asu\be = sv$ and then $\ova\ove\ovbe = \ove$, where $\ovbe = \be + cR$. Similarly, $\ovb\ovx(\ov1 - \ove) = \ov1 - \ove$ for some element $\ovx \in \ovR$. We have proved that if $\ova\ovR + \ovb\ovR = \ovR$, then there exists an idempotent $\ove$ such that $\ove \in \ova\ovR$ and $\ov1 - \ove \in \ovb\ovR$ i. e. $R$ is an exchange ring.
\end{Proof}

\begin{prop}\label{p9}
      Let $R$ be a commutative Bezout ring of diadic range 1. Then for any divisor $\al$ of a diadem $a + b\la$ and elements $c, d \in R$ such that $\al R + cR + dR = R$ there exists an element $\mu \in R$ such that $\al R + (c + d\mu)R = R$.
\end{prop}

\begin{Proof}
      Let $a + b\la = \al\be$ and $cR + dR = kR$. Since $R$ is a ring of diadic range one, by Theorem \ref{t1} and Theorem \ref{t3} we have that $c = kc_1$, $d = kd_1$ for some elements $c_1, d_1 \in R$ such that $c_1R + d_1R = R$.

      Since $(a + b\la)R + c_1R + d_1R = R$ there exists an element $\mu \in R$ such that $(a + b\la)R + (c_1 + d_1\mu)R = R$. Whence $\al R + (c_1 + d_1\mu)R = R$. Since $\al R + cR + dR = R$ and $cR + dR = kR$, we have $\al R + kR = R$. The equality $\al R + (c_1 + d_1\mu)R = R$ implies $R = \al R + k(c_1 + d_1\mu)R = \al R + (c + d\mu)R$.
\end{Proof}

\begin{thm}\label{t10}
      A commutative Bezout ring is an elementary divisor ring if and and only if it is a ring of dyadic range 1.
\end{thm}

\begin{Proof}
      Let $R$ be a commutative Bezout ring of dyadic range 1. By Theorem \ref{t1} and Theorem \ref{t3} we have that $R$ is an Hermite ring. Suppose that
      $A=\begin{pmatrix}
         a & 0 \\
         b & c
      \end{pmatrix} \in M_2(R)$ with $aR + bR + cR = R$. By \cite{KaplEDMod} it suffices to check that $A$ admits an elementary reduction.

      Since $R$ is a ring of dyadic range 1, there exist some elements $x, y \in R$ such that $b + ax + cy = w$ is diadem. Then
      $\begin{pmatrix}
         1 & 0 \\
         x & 1
      \end{pmatrix}
      \begin{pmatrix}
         a & 0 \\
         b & c
      \end{pmatrix}
      \begin{pmatrix}
         1 & 0 \\
         y & 1
      \end{pmatrix} =
      \begin{pmatrix}
         a & 0 \\
         w & c
      \end{pmatrix}.$
      Obviously, $aR + wR + cR = R$. Since $R$ is an Hermite ring, there exists an invertible matrix $Q \in GL_2(R)$ such that $(w, c)Q = (\al, 0)$ for some element $\al \in R$. Obviously, $\al$ is a divisor of the diadem $w$.
      Let
      $\begin{pmatrix}
         a & 0 \\
         w & c
      \end{pmatrix}Q =
      \begin{pmatrix}
         a' & c' \\
         0 & \al
      \end{pmatrix}$
      with some elements $a', c' \in R$. It is easily seen that $\al R + a'R + c'R = R$. By Proposition \ref{p9}, there exists an element $\mu \in R$ such that $\al R + (c' + a'\mu)R = R$. Then
      $\begin{pmatrix}
         a' & c' \\
         0 & \al
      \end{pmatrix}
      \begin{pmatrix}
         \mu & 1 \\
         1 & 0
      \end{pmatrix} =
      \begin{pmatrix}
         c' + a'\mu & a' \\
         \al & 0
      \end{pmatrix}$.
      Since $R$ is an Hermite ring and $\al R + (c' + a'\mu)R = R$, there exist an invertible matrix $P \in GL_2(R)$ such that
      $P\begin{pmatrix}
         c' + a'\mu & a' \\
         \al & 0
      \end{pmatrix} =
      \begin{pmatrix}
         1 & \ast \\
         \ast & \ast
      \end{pmatrix} = B$.
      Obviously, the matrix $B$ admits a diagonal reduction. Therefore,
      $A = \begin{pmatrix}
         a & 0 \\
         b & c
      \end{pmatrix}$
      admits a diagonal reduction and therefore $R$ is an elementary divisor ring.

      Let $R$ be an elementary divisor ring. By \cite{Zab22}, for some elements $a, b \in R$ such that $aR + bR = R$ there exists an element $\la \in R$ such that for any element $c \in R$ we have $a + b\la = uv$, where $uR + cR = R$, $vR + (1 - c)R = R$ and $uR + vR = R$. By Proposition \ref{p7}, $\ovR = R/(a + b\la)R$ is an exchange ring. Since $R$ is a commutative ring, by \cite{Ch9}, $R$ is of idempotent stable range 1, i.e. stable range $\ovR$ is equal 1. By Proposition \ref{p6}, the element $a + b\la$ is a diadem.
\end{Proof}

As a consequence of this theorem we obtain the following result.

\begin{prop}\label{p11}
      Let $R$ be a commutative Bezout ring of diadic range 1. Then for any ideal $I$ of $R$ the factor-ring $R/I$ is a ring of diadic range 1.
\end{prop}

\begin{Proof}
      Since any homomorphic image of an elementary divisor ring is an elementary divisor ring, by Theorem \ref{t10} we have proved our proposition.
\end{Proof}

Moreover, we obtain the result

\begin{prop}\label{p12}
      Let $R$ be a commutative semi hereditary Bezout ring. If for any regular element (non zero divisor) $r \in R$ the factor-ring $R/rR$ is a ring of dyadic range 1 then $R$ is a ring of dyadic range 1.
\end{prop}

\begin{Proof}
      By \cite{Shor1974}, if for a commutative semihereditary Bezout ring $R$ the factor-ring $R/rR$ for any regular element $r \in R$ is an elementary divisor ring, then $R$ is an elementary divisor ring. By Theorem \ref{t10}, Proposition \ref{p12} is obvious
\end{Proof}

Consequently, we obtain that an example of a commutative Hermite ring which is an elementary divisor ring \cite{Gill} is an example of a commutative Bezout ring of stable range 2 which is not a ring of diadic range 1.

\vspace{5ex}
\renewcommand\refname{References}

\end{document}